\documentclass[11pt,a4paper]{amsart}
\usepackage[T1]{fontenc}
\usepackage{lmodern,microtype}
\usepackage[margin=27mm]{geometry}
\usepackage{amsmath,amssymb,mathtools}
\usepackage{setspace}
\usepackage[hidelinks]{hyperref}

\newtheorem{lemma}{Lemma}
\newtheorem{theorem}[lemma]{Theorem}
\newtheorem{proposition}[lemma]{Proposition}
\newtheorem{corollary}[lemma]{Corollary}
\newcommand{\HH}{\mathbb H}
\newcommand{\RR}{\mathbb R}
\newcommand{\CC}{\mathbb C}
\newcommand{\Rea}{\operatorname{Re}}
\newcommand{\Ima}{\operatorname{Im}}
\newcommand{\Bplus}{B^+}
\newcommand{\Cplus}{C^+}
\newcommand{\norm}[1]{\lVert#1\rVert}
\DeclareMathOperator{\diag}{diag}
\DeclareMathOperator{\spann}{span}
\title{On the convexity of the upper bild}
\author{Lu\'{\i}s Carvalho}
\address{ISCTE -- Lisbon University Institute, Av. das For\c{c}as Armadas,
1649-026 Lisbon, Portugal; and Center for Mathematical Analysis,
Geometry and Dynamical Systems, Department of Mathematics,
Instituto Superior T\'{e}cnico, Universidade de Lisboa,
Av. Rovisco Pais 1, 1049-001 Lisbon, Portugal}
\email{luis.carvalho@iscte-iul.pt}
\date{}

\begin{document}
\setstretch{1.25}
\begin{abstract}
We give a proof of the convexity of the upper bild of a quaternionic
matrix, a theorem of So and Thompson. After reduction to $2\times2$
matrices, we introduce real coordinates on a sphere and obtain two
representations of the upper bild over the closed unit ball,
involving a convex function and a concave function. Minimizing
the first and maximizing the second over
parameters with a fixed real part yield a convex lower boundary
function and a concave upper boundary function, from which
convexity follows.
\end{abstract}
\maketitle

\section{Introduction}\label{sec:introduction}

Let $\HH$ denote the skew field of Hamilton quaternions and let
$M_n(\HH)$ be the set of $n\times n$ matrices with quaternionic entries.
Vectors are columns, with scalars acting on the right.
For $A\in M_n(\HH)$, its \emph{numerical range} is
\[
 W_{\HH}(A)=\{\mathbf x^*A\mathbf x:
              \mathbf x\in\HH^n,\ \norm{\mathbf x}=1\},
 \qquad \norm{\mathbf x}^2=\mathbf x^*\mathbf x.
\]
In the complex setting, the Toeplitz--Hausdorff theorem guarantees
that the numerical range is convex; see, for instance, \cite{GR}.
For quaternionic matrices convexity may fail, but the numerical range
has a symmetry that allows its geometry to be studied in a complex
plane. To describe this planar set, fix an imaginary unit $i$,
identify $\CC$ with $\RR+\RR i$, and define the \emph{upper bild} by
\[
 \Bplus(A)=W_{\HH}(A)\cap\CC^+,
 \qquad \CC^+=\{a+ib:a\in\RR,\ b\geq0\}.
\]
Conjugating any quaternionic numerical value by a unit quaternion
gives another value in the numerical range. Each such similarity
class has exactly one representative in $\CC^+$, determined by its
real part and the norm of its imaginary part. Consequently, the
upper bild determines the whole numerical range: one recovers
$W_{\HH}(A)$ by taking all quaternions similar to the elements of
the upper bild.

So and Thompson \cite{ST} proved that the upper bild is convex for
every quaternionic matrix. Their theorem provides a convex planar
set even when the numerical range itself is not convex, and makes
the geometry of the upper bild central to the study of quaternionic
numerical ranges. 

In this paper we give a proof of the So--Thompson theorem based on
real coordinates and elementary properties of convex functions.
Our aim is to show how the formulas for the numerical values yield
a convex lower boundary function and a concave upper boundary
function. We prove the following established result.

\begin{theorem}\label{thm:upper-bild-convexity}
For every $A\in M_n(\HH)$, the upper bild $\Bplus(A)$ is convex.
\end{theorem}

The proof begins with the classical reduction to $2\times2$ matrices.
After choosing an orthonormal basis in which the skew-Hermitian part
is diagonal, we introduce real coordinates that describe the upper
bild as the set of complex values obtained from points on the unit
sphere in $\RR^5$.

We show that the same formula can be evaluated throughout the closed
unit ball without introducing any new values: every complex number
obtained from a point inside the ball is already obtained from some
point on the sphere. Thus, extending the domain from the sphere to
the ball leaves the set of values unchanged.

We then find a second expression for the imaginary part, giving
another representation of the upper bild over the same ball. Both
representations produce exactly the same set of complex values, but
the first expression for the imaginary part is convex and the second
is concave. These two properties allow us to study the lower and
upper boundaries of the upper bild and establish its convexity.

\medskip
\noindent\textit{Use of large language models.}
A large language model was used to support mathematical discussion,
the development and checking of arguments and calculations,
literature searches, and the drafting and revision of the manuscript.
Its use therefore included mathematical reasoning and exposition
as well as language editing. The author takes responsibility for
the mathematical claims, references and final text.

\section{Preliminaries}\label{sec:preliminaries}

We recall the quaternionic facts used in the proof; for background
on quaternionic matrices, see \cite{Zh}.
Two quaternions $p,q$ are \emph{equivalent}, or \emph{similar}, if
$p=\omega^*q\omega$ for some unit quaternion $\omega$. The equivalence class of $q$ is
\[
 [q]=\{\omega^*q\omega:\omega\in\HH,\ |\omega|=1\}.
\]
Conjugation by a unit quaternion fixes the real part and rotates the
imaginary part, preserving its norm:
\[
 \Rea(\omega^*q\omega)=\Rea q,
 \qquad
 \Ima(\omega^*q\omega)=\omega^*(\Ima q)\omega.
\]
Every direction of the imaginary part can be obtained in this way.
Consequently, $p$ and $q$ are equivalent if and only if
\[
 \Rea p=\Rea q,
 \qquad |\Ima p|=|\Ima q|.
\]
Thus $[q]$ is a sphere with centre $\Rea q$ and radius $|\Ima q|$
in the affine space $\Rea q+\Ima\HH$. For real $q$, it reduces to
the singleton $\{q\}$. Each class has exactly one representative in
$\CC^+$, namely
\[
 [q]\cap\CC^+=\{\Rea q+i|\Ima q|\}.
\]

Let $z=\mathbf x^*A\mathbf x$ with $\norm{\mathbf x}=1$.
For every unit quaternion $\omega$, the vector $\mathbf x\omega$
also has norm one and
\[
 (\mathbf x\omega)^*A(\mathbf x\omega)
 =\omega^*z\omega\in[z].
\]
As $\omega$ ranges over the unit quaternions, $\omega^*z\omega$
ranges over the whole equivalence class $[z]$. Hence
$[z]\subseteq W_{\HH}(A)$, and its representative
$\Rea z+i|\Ima z|$ belongs to $\Bplus(A)$.

Convexity of the upper bild reduces to the
$2\times2$ case by the classical Toeplitz--Hausdorff compression argument.

\begin{lemma}\label{lem:compression}
If $\Bplus(C)$ is convex for every $C\in M_2(\HH)$, then
$\Bplus(A)$ is convex for every $A\in M_n(\HH)$.
\end{lemma}
\begin{proof}
Take $z_1,z_2\in\Bplus(A)$ and choose $\mathbf u,\mathbf v$ with
$\norm{\mathbf u}=\norm{\mathbf v}=1$ such that
$z_1=\mathbf u^*A\mathbf u$ and $z_2=\mathbf v^*A\mathbf v$. Put
$E=\spann_{\HH}\{\mathbf u,\mathbf v\}$.
Let $P_E$ be the orthogonal projection onto $E$ and let $C=P_EA|_E$.
For $\mathbf x\in E$,
\[
 \mathbf x^*C\mathbf x=\mathbf x^*A\mathbf x.
\]
Hence
\[
 z_1,z_2\in\Bplus(C)\subseteq\Bplus(A).
\]
If $\dim_{\HH}E=2$, the hypothesis gives
$[z_1,z_2]\subseteq\Bplus(C)$.
If $\dim_{\HH}E=1$, all unit vectors in $E$ differ by a unit right
scalar, so $\Bplus(C)$ is a singleton and $z_1=z_2$.
In either case, $[z_1,z_2]\subseteq\Bplus(A)$.
\end{proof}

By Lemma~\ref{lem:compression}, proving convexity for $2\times2$
matrices will establish it in every dimension, and we can simplify
the matrices under consideration by choosing an orthonormal basis
in which their skew-Hermitian part is diagonal. Such a change of
basis preserves the numerical range and the upper bild because,
for every unitary matrix $U\in M_2(\HH)$, we have
\[
 \mathbf x^*(U^*AU)\mathbf x=(U\mathbf x)^*A(U\mathbf x),
\]
and $U$ maps the unit sphere onto itself. Writing
$\mathbb S_{\HH^2}=\{\mathbf x\in\HH^2:\norm{\mathbf x}=1\}$,
this means that
\[
 U\mathbb S_{\HH^2}
 =U\{\mathbf x\in\HH^2:\norm{\mathbf x}=1\}
 =\mathbb S_{\HH^2}.
\]
Consequently,
\[
 \Bplus(A)=\Bplus(U^*AU).
\]

To choose this basis for $A\in M_2(\HH)$, we first separate its
Hermitian and skew-Hermitian parts by writing
\[
 A=H_0+S_0,\qquad H_0=\frac{A+A^*}{2},
 \qquad S_0=\frac{A-A^*}{2}.
\]
Since $S_0^*=-S_0$, the quaternionic spectral theorem
\cite[Corollary~6.2]{Zh} allows us to
diagonalize $S_0$ in an orthonormal basis, with purely imaginary
complex diagonal entries chosen to have nonnegative imaginary
parts; thus there is a unitary matrix $U$ such that
\[
 U^*S_0U=\diag(s_1i,s_2i)=Si,
 \qquad s_1,s_2\geq0.
\]
In this same basis, the matrix $H=U^*H_0U$ remains Hermitian, so
its diagonal entries are real and its off-diagonal entries are
conjugates of one another. Replacing $A$ by $U^*AU$, which leaves
the upper bild unchanged, we may therefore work with the form
\[
 A=H+Si,\qquad
 H=\begin{pmatrix}h_1&w\\\bar w&h_2\end{pmatrix},
 \qquad
 S=\begin{pmatrix}s_1&0\\0&s_2\end{pmatrix},
\]
where $h_1,h_2\in\RR$, $w\in\HH$, and $s_1,s_2\geq0$.
We shall use this representation throughout, with $H$ denoting
the Hermitian part and the real diagonal matrix $S$ recording the
coefficients of the skew-Hermitian part $Si=iS$.

Multiplying a vector on the right by a unit quaternion allows us to
choose a real and nonnegative first entry while preserving its
upper-bild representative. We collect the unit vectors with this
property in the set
\[
 \Sigma=\{\mathbf x=(r,q)^T\in\HH^2:
                  r\in[0,1],\ r^2+|q|^2=1\}.
\]

\begin{lemma}\label{lem:real-first-entry}
For $A=H+Si$ as above,
\[
 \Bplus(A)=
 \left\{\mathbf x^*H\mathbf x+i\bigl|\mathbf x^*(Si)\mathbf x\bigr|:
 \mathbf x\in\Sigma\right\}.
\]
\end{lemma}
\begin{proof}
Let $\mathbf x=(x_1,x_2)^T\in\HH^2$ be any unit vector.
Since $H$ is Hermitian and $Si$ is skew-Hermitian,
$\mathbf x^*H\mathbf x$ is real and $\mathbf x^*(Si)\mathbf x$ is
purely imaginary. The upper-bild representative of
$\mathbf x^*A\mathbf x$ is therefore
\[
 \mathbf x^*H\mathbf x+i\bigl|\mathbf x^*(Si)\mathbf x\bigr|.
\]
If $x_1\neq0$, set
\[
 \omega=\frac{\overline{x_1}}{|x_1|},\qquad
 \mathbf y=\mathbf x\omega.
\]
Then $|\omega|=1$, $\norm{\mathbf y}=1$, and $y_1=|x_1|$,
so $\mathbf y\in\Sigma$.
Moreover,
\[
 \mathbf y^*A\mathbf y
 =\bar\omega(\mathbf x^*A\mathbf x)\omega.
\]
Quaternionic similarity preserves the real part and the norm of the
imaginary part, so $\mathbf y$ gives the same upper-bild representative
as $\mathbf x$. If $x_1=0$, we already have $\mathbf x\in\Sigma$
and may take $\mathbf y=\mathbf x$.
\end{proof}

\section{Convexity of the upper bild}\label{sec:convexity}

To choose coordinates adapted to the numerical values, we first
examine the quadratic form. For $\mathbf x=(r,q)^T\in\Sigma$,
expanding $A=H+Si$ gives
\begin{equation}\label{eq:quadratic-decomposition}
 \mathbf x^*A\mathbf x
 =\underbrace{h_1r^2+h_2|q|^2+2r\Rea(wq)}_{\mathbf x^*H\mathbf x}
  +\underbrace{s_1r^2i+s_2\bar q i q}_{\mathbf x^*(Si)\mathbf x}.
\end{equation}
The first group is real and the second is purely imaginary, so
the real part depends on the two squared norms $r^2,|q|^2$ and
the product $rq$. Since $r$ is real, the mixed term satisfies
\[
 2r\Rea(wq)=\Rea\bigl(w(2rq)\bigr),
\]
which suggests retaining the quaternion $2rq$.

The condition $r^2+|q|^2=1$ means that the two squared norms
are determined by their difference $t=r^2-|q|^2$, since adding
and subtracting these two relations gives
\begin{equation}\label{eq:squared-norms-from-t}
 r^2=\frac{1+t}{2},\qquad |q|^2=\frac{1-t}{2}.
\end{equation}
Substituting into the real part of~\eqref{eq:quadratic-decomposition}
therefore yields
\begin{equation}\label{eq:real-part-before-phi}
 \Rea(\mathbf x^*A\mathbf x)
 =\frac{h_1+h_2}{2}+\frac{h_1-h_2}{2}\,t
  +\Rea\bigl(w(2rq)\bigr).
\end{equation}
This expression is affine in $t$ and in the four real coefficients
of $2rq$. These quantities thus provide five real coordinates
in which the real part takes a particularly simple form.

Since $t=r^2-|q|^2$ and $r^2+|q|^2=1$, we have
\[
 t^2+|2rq|^2
 =(r^2-|q|^2)^2+4r^2|q|^2
 =(r^2+|q|^2)^2=1.
\]
Thus these five real quantities lie on the unit sphere in $\RR^5$.
We now introduce the coordinate map precisely and show that every
point of this sphere is attained.

To write the quaternion $2rq$ in real coordinates, choose an imaginary unit $j$
orthogonal to $i$, put $k=ij$, and define
\[
 \psi:\HH\longrightarrow\RR^4,\qquad
 \psi(p_0+p_1i+p_2j+p_3k)=(p_0,p_1,p_2,p_3).
\]
The map $\psi$ is a real-linear bijection and an isometry, since
the quaternionic norm and the Euclidean norm are both the square
root of the sum of squares of these four real coefficients.
For $\mathbf p=(p_0,p_1,p_2,p_3)\in\RR^4$, its inverse therefore
satisfies
\[
 \psi^{-1}(\mathbf p)=p_0+p_1i+p_2j+p_3k,
 \qquad |\psi^{-1}(\mathbf p)|=\norm{\mathbf p}.
\]
For $\mathbf x=(r,q)^T\in\Sigma$, we therefore define
$\phi:\Sigma\longrightarrow\RR^5$ by
\begin{equation}\label{eq:phi}
 \phi(\mathbf x)=\bigl(r^2-|q|^2,\,\psi(2rq)\bigr).
\end{equation}
We write $\boldsymbol\xi=\phi(\mathbf x)=(t,\mathbf p)$, where
\[
 t=r^2-|q|^2,\qquad
 \mathbf p=(p_0,p_1,p_2,p_3)=\psi(2rq).
\]
Thus $t$ records the difference of the two squared norms, while
$\mathbf p$ records the four real coefficients of $2rq$; whenever
this quaternion is needed in a product, it is recovered as
$\psi^{-1}(\mathbf p)=2rq$.
The map is continuous because its components are real quadratic
polynomials in $r$ and the real coordinates of $q$. Writing
$S^4=\{\boldsymbol\xi\in\RR^5:\norm{\boldsymbol\xi}=1\}$,
we obtain the following description of its image.

\begin{lemma}\label{lem:phi}
\[
 \phi(\Sigma)=S^4.
\]
\end{lemma}
\begin{proof}
For every $\mathbf x=(r,q)^T\in\Sigma$, the isometry property of
$\psi$ gives
\[
 \begin{aligned}
 \norm{\phi(\mathbf x)}^2
 &=\norm{\bigl(r^2-|q|^2,\,\psi(2rq)\bigr)}^2\\
 &=(r^2-|q|^2)^2+4r^2|q|^2
 =(r^2+|q|^2)^2=1,
 \end{aligned}
\]
which proves $\phi(\Sigma)\subseteq S^4$.
Conversely, given
$\boldsymbol\xi=(t,\mathbf p)\in S^4$, where
$\mathbf p=(p_0,p_1,p_2,p_3)\in\RR^4$,
the same property gives
$|\psi^{-1}(\mathbf p)|^2=\norm{\mathbf p}^2=1-t^2$.
If $t>-1$, choose
\[
 r=\sqrt{\frac{1+t}{2}},\qquad
 q=\frac{p_0+p_1i+p_2j+p_3k}{2r}.
\]
Then $|q|^2=(1-t)/2$, which gives $r^2+|q|^2=1$,
$r^2-|q|^2=t$, and $\psi(2rq)=\mathbf p$. Hence, with
$\mathbf x=(r,q)^T$, we have $\mathbf x\in\Sigma$ and
$\phi(\mathbf x)=\boldsymbol\xi$.
If $t=-1$, we have $\mathbf p=\mathbf 0$, and the same conclusion follows by
choosing $\mathbf x=(0,1)^T\in\Sigma$.
\end{proof}

By Lemma~\ref{lem:real-first-entry}, it suffices to calculate the
real part and the norm of the imaginary part of
$\mathbf x^*A\mathbf x$ for $\mathbf x\in\Sigma$, and we now
express both quantities in terms of the five real coordinates
$\boldsymbol\xi=(t,\mathbf p)=\phi(\mathbf x)$.

For $\boldsymbol\xi=(t,\mathbf p)\in\RR^5$, with
$\mathbf p=(p_0,p_1,p_2,p_3)$, define
\begin{align}
 h(\boldsymbol\xi)
 &=\frac{h_1+h_2}{2}+\frac{h_1-h_2}{2}\,t
   +\Rea\bigl(w\psi^{-1}(\mathbf p)\bigr),
 \label{eq:h-coordinates}\\
 \rho(\boldsymbol\xi)
 &=\sqrt{\left(\frac{s_1(1+t)-s_2(1-t)}{2}\right)^2
                   +s_1s_2(p_0^2+p_1^2)}.
 \label{eq:rho-coordinates}
\end{align}
Since $\psi^{-1}$ is real-linear, so is the map
$\mathbf p\mapsto\Rea\bigl(w\psi^{-1}(\mathbf p)\bigr)$,
which shows that $h$ is affine in the five real coordinates.
The expression defining $\rho$ is nonnegative because
$s_1s_2\geq0$, so $\rho$ is a well-defined nonnegative function.

\begin{lemma}\label{lem:coordinate-values}
For every $\mathbf x\in\Sigma$,
\[
 \Rea(\mathbf x^*A\mathbf x)=h\bigl(\phi(\mathbf x)\bigr),
 \qquad
 |\Ima(\mathbf x^*A\mathbf x)|=\rho\bigl(\phi(\mathbf x)\bigr).
\]
\end{lemma}
\begin{proof}
Write $\mathbf x=(r,q)^T$ and
$\boldsymbol\xi=\phi(\mathbf x)=(t,\mathbf p)$.
Since $\psi^{-1}(\mathbf p)=2rq$, the expression for the real
part in~\eqref{eq:real-part-before-phi} agrees with the definition
of $h$ in~\eqref{eq:h-coordinates}, proving the first identity.

For the second identity, the decomposition
in~\eqref{eq:quadratic-decomposition} gives
\[
 \begin{aligned}
 |\Ima(\mathbf x^*A\mathbf x)|^2
 &=\bigl|s_1r^2i+s_2\bar q i q\bigr|^2\\
 &=s_1^2r^4+s_2^2|q|^4
   +2s_1s_2r^2\Rea(-i\bar q i q).
 \end{aligned}
\]
To express $4r^2\Rea(-i\bar q i q)$ in terms of the same four
real coordinates, write
\[
 2rq=p_0+p_1i+p_2j+p_3k=u+vj,
 \qquad u=p_0+p_1i,\quad v=p_2+p_3i\in\CC.
\]
For complex coefficients, $ju=\bar u j$ and $jv=\bar vj$;
together with $ji=-ij$ and $j^2=-1$, these identities give
\[
 -i\,\overline{u+vj}\,i
 =-i(\bar u-vj)i=\bar u+vj
\]
and
\[
 \begin{aligned}
 (\bar u+vj)(u+vj)
 &=|u|^2+\bar u vj+v\bar u j-v\bar v\\
 &=|u|^2-|v|^2+2\bar u vj.
 \end{aligned}
\]
Since $\bar u v\in\CC$, the term $2\bar u vj$ is a real linear
combination of $j$ and $k$ and therefore has zero real part.
Using that $r$ is real, we consequently obtain
\[
 \begin{aligned}
 4r^2\Rea(-i\bar q i q)
 &=\Rea\bigl(-i\,\overline{2rq}\,i\,(2rq)\bigr)\\
 &=\Rea\bigl((\bar u+vj)(u+vj)\bigr)\\
 &=|u|^2-|v|^2\\
 &=p_0^2+p_1^2-p_2^2-p_3^2.
 \end{aligned}
\]
Substituting this identity and~\eqref{eq:squared-norms-from-t}
into the squared norm gives
\[
 |\Ima(\mathbf x^*A\mathbf x)|^2
 =\frac{s_1^2(1+t)^2+s_2^2(1-t)^2}{4}
  +\frac{s_1s_2}{2}(p_0^2+p_1^2-p_2^2-p_3^2).
\]
Finally, $\boldsymbol\xi\in S^4$ gives
$p_2^2+p_3^2=1-t^2-p_0^2-p_1^2$, so the preceding expression becomes
\[
 \begin{aligned}
 |\Ima(\mathbf x^*A\mathbf x)|^2
 &=\left(\frac{s_1(1+t)-s_2(1-t)}{2}\right)^2
   +s_1s_2(p_0^2+p_1^2)\\
 &=\rho\bigl(\phi(\mathbf x)\bigr)^2.
 \end{aligned}
\]
Taking nonnegative square roots proves the second identity.

When $r=0$, the condition $\mathbf x\in\Sigma$ gives $|q|=1$,
so $t=-1$ and $\mathbf p=\mathbf0$. Substituting these coordinates
into the definitions of $h$ and $\rho$ yields
\[
 h(\boldsymbol\xi)=h_2,\qquad
 \rho(\boldsymbol\xi)=\sqrt{s_2^2}=s_2,
\]
where the last equality uses $s_2\geq0$. Direct evaluation of the
original quadratic form gives
\[
 \mathbf x^*A\mathbf x=h_2+s_2\bar q i q.
\]
Since $\bar q i q$ is purely imaginary and has norm
$|q|^2=1$, this value has real part $h_2$ and imaginary part of
norm $s_2$, in agreement with the coordinate formulas.
\end{proof}

The preceding calculation leads to a set defined entirely
in terms of the five real coordinates and the parameters of the matrix $A$:
\begin{equation}\label{eq:Cplus}
 \Cplus(A)=\{h(\boldsymbol\xi)+i\rho(\boldsymbol\xi):
                      \boldsymbol\xi\in S^4\},
\end{equation}
where $h$ and $\rho$ are given by
\eqref{eq:h-coordinates} and~\eqref{eq:rho-coordinates}.

\begin{proposition}\label{thm:coordinate-model}
For $A=H+Si$ in the form fixed above,
\[
 \Bplus(A)=\Cplus(A).
\]
\end{proposition}
\begin{proof}
For $z\in\Bplus(A)$, Lemma~\ref{lem:real-first-entry} supplies
$\mathbf x\in\Sigma$ such that $z$ is the upper-bild representative
of $\mathbf x^*A\mathbf x$. Lemma~\ref{lem:coordinate-values}
and the inclusion $\phi(\Sigma)\subseteq S^4$ then give
\[
 z=h(\phi(\mathbf x))+i\rho(\phi(\mathbf x))\in\Cplus(A).
\]
Conversely, given $\boldsymbol\xi\in S^4$,
Lemma~\ref{lem:phi} supplies $\mathbf x\in\Sigma$ with
$\phi(\mathbf x)=\boldsymbol\xi$, and
Lemma~\ref{lem:coordinate-values} identifies
$h(\boldsymbol\xi)+i\rho(\boldsymbol\xi)$ as the unique element
of $[\mathbf x^*A\mathbf x]$ in $\CC^+$. Since the numerical range
contains this whole equivalence class, that element belongs to
$\Bplus(A)$.
\end{proof}

Having identified the upper bild with $\Cplus(A)$, consider the map
\[
 \gamma:S^4\longrightarrow\CC^+,\qquad
 \gamma(\boldsymbol\xi)=h(\boldsymbol\xi)+i\rho(\boldsymbol\xi).
\]
Equation~\eqref{eq:Cplus} says that $\Cplus(A)=\gamma(S^4)$.
We now extend $\gamma$ to the closed unit ball
\[
 D^5=\{\boldsymbol\xi\in\RR^5:\norm{\boldsymbol\xi}\leq1\}.
\]
The same formula defines this extension, since $h$ and $\rho$
are already defined on $\RR^5$. Our aim is to prove that
$\gamma(D^5)=\gamma(S^4)$: enlarging the domain from the sphere
to the ball does not enlarge the set of values. Working on the
ball will then allow us to take convex combinations while remaining
in the domain. We begin with the following fact about affine maps.

\begin{lemma}\label{lem:affine-sphere-ball}
For every affine map $L:\RR^5\to\RR^4$,
\[
 L(S^4)=L(D^5).
\]
\end{lemma}
\begin{proof}
Write $L(\mathbf u)=T\mathbf u+\mathbf b$, where $T$ is linear.
Since $T$ maps $\RR^5$ into $\RR^4$, its kernel contains a
nonzero vector $\mathbf v$. For any $\mathbf u\in D^5$, we
therefore have
\[
 L(\mathbf u+\lambda\mathbf v)=L(\mathbf u)
 \qquad\text{for every }\lambda\in\RR.
\]
At $\lambda=0$, we have
$\norm{\mathbf u+\lambda\mathbf v}=\norm{\mathbf u}\leq1$
because $\mathbf u\in D^5$. For $\lambda\geq0$, the triangle inequality gives
\[
 \lambda\norm{\mathbf v}
 =\norm{(\mathbf u+\lambda\mathbf v)-\mathbf u}
 \leq\norm{\mathbf u+\lambda\mathbf v}+\norm{\mathbf u},
\]
and hence
\[
 \norm{\mathbf u+\lambda\mathbf v}
 \geq\lambda\norm{\mathbf v}-\norm{\mathbf u}
 \longrightarrow\infty
 \qquad\text{as }\lambda\to\infty,
\]
since $\mathbf v\ne\mathbf0$. By continuity, there is
$\lambda_0\geq0$ such that
$\boldsymbol\eta=\mathbf u+\lambda_0\mathbf v\in S^4$,
and $L(\boldsymbol\eta)=L(\mathbf u)$. Thus every value
attained on $D^5$ is attained on $S^4$, and the reverse inclusion
follows from $S^4\subseteq D^5$.
\end{proof}

To apply this lemma, write a general point of $\RR^5$ as
$\boldsymbol\xi=(t,p_0,p_1,p_2,p_3)$. The value
$\gamma(\boldsymbol\xi)$ for $\boldsymbol\xi\in D^5$ is
determined by the real number $h(\boldsymbol\xi)$ and the first
three coordinates $t,p_0,p_1$, which determine
$\rho(\boldsymbol\xi)$ through~\eqref{eq:rho-coordinates}.
Geometrically, this leaves $p_2,p_3$ free to vary without changing
$\rho$. Once $t,p_0,p_1$ are fixed, keeping $h$ unchanged imposes
at most one real linear condition on variations of $p_2,p_3$,
because $h$ is affine. There is therefore a nonzero direction
along which both $h$ and $\rho$ remain constant. The line through
any point of $D^5$ in this direction meets $S^4$, so the same
complex value is attained on the sphere.

We therefore express $\gamma$ as a composition: an affine map
first records $h(\boldsymbol\xi),t,p_0,p_1$, and a second map
uses these four real numbers to calculate the complex value.
The lemma gives the same set of values for the first map on
the sphere and on the ball; applying the second map preserves
this equality.

\begin{proposition}\label{thm:ball-representations}
For the fixed matrix $A=H+Si$,
\begin{equation}\label{eq:ball-representations}
 \Cplus(A)=\gamma(D^5)=\{h(\boldsymbol\xi)+i\rho(\boldsymbol\xi):
                                  \boldsymbol\xi\in D^5\}.
\end{equation}
\end{proposition}
\begin{proof}
For $\boldsymbol\xi=(t,p_0,p_1,p_2,p_3)\in\RR^5$, recall that
\[
 h(\boldsymbol\xi)
 =\frac{h_1+h_2}{2}+\frac{h_1-h_2}{2}\,t
   +\Rea\bigl(w(p_0+p_1i+p_2j+p_3k)\bigr).
\]
Since multiplication by $w$ and taking real parts are real-linear,
$h$ is affine in $\boldsymbol\xi$. We therefore define the affine map
\[
 L:\RR^5\longrightarrow\RR^4,\qquad
 L(\boldsymbol\xi)=\bigl(h(\boldsymbol\xi),t,p_0,p_1\bigr).
\]
Its remaining components are coordinate projections, which are linear.
Next define $G:\RR^4\to\CC^+$ by
\[
 G(h,t,p_0,p_1)
 =h+i\sqrt{\left(\frac{s_1(1+t)-s_2(1-t)}{2}\right)^2
                  +s_1s_2(p_0^2+p_1^2)}.
\]
Its square root is real because $s_1s_2\geq0$, and substituting
$L(\boldsymbol\xi)$ into this formula gives
\[
 (G\circ L)(\boldsymbol\xi)
 =h(\boldsymbol\xi)+i\rho(\boldsymbol\xi)=\gamma(\boldsymbol\xi)
 \qquad\text{for every }\boldsymbol\xi\in D^5.
\]
Since Lemma~\ref{lem:affine-sphere-ball} gives $L(S^4)=L(D^5)$,
applying $G$ to these equal sets and using~\eqref{eq:Cplus} yields
\[
 \begin{aligned}
 \Cplus(A)
 &=\gamma(S^4)=G\bigl(L(S^4)\bigr)\\
 &=G\bigl(L(D^5)\bigr)=\gamma(D^5)\\
 &=\{h(\boldsymbol\xi)+i\rho(\boldsymbol\xi):
                                      \boldsymbol\xi\in D^5\}.
 \end{aligned}
\]
\end{proof}

To study the lower and upper boundaries of $\Cplus(A)$, we will
fix the real part and consider the smallest and largest imaginary
parts that occur. As we shall show, $\rho$ is convex on $D^5$,
and minimizing it over parameters with a fixed real part gives
a convex lower boundary function. For the upper boundary, we
seek a second radius function that is concave on the ball and
that, together with the same real part $h$, gives the same set
of complex values. Maximizing this function over parameters
with a fixed real part will then give a concave upper boundary
function.

We return to the sphere to obtain a second expression for the
radius. For $\boldsymbol\xi=(t,p_0,p_1,p_2,p_3)\in S^4$,
recall from~\eqref{eq:rho-coordinates} that
\[
 \rho(\boldsymbol\xi)
 =\sqrt{\left(\frac{s_1(1+t)-s_2(1-t)}{2}\right)^2
                   +s_1s_2(p_0^2+p_1^2)}.
\]
The condition $\norm{\boldsymbol\xi}=1$ means that
\[
 t^2+p_0^2+p_1^2+p_2^2+p_3^2=1,
\]
so we can replace the sum $p_0^2+p_1^2$ using
\[
 p_0^2+p_1^2=1-t^2-p_2^2-p_3^2.
\]
Substituting gives
\[
 \begin{aligned}
 \rho(\boldsymbol\xi)^2
 &=\left(\frac{s_1(1+t)-s_2(1-t)}{2}\right)^2
       +s_1s_2(1-t^2)-s_1s_2(p_2^2+p_3^2)\\
 &=\left(\frac{s_1(1+t)+s_2(1-t)}{2}\right)^2
       -s_1s_2(p_2^2+p_3^2).
 \end{aligned}
\]
The last equality follows by expanding the squares and simplifying.
Taking the nonnegative square root therefore gives another
expression for the radius on $S^4$, which we denote by
\begin{equation}\label{eq:rho-tilde}
 \widetilde\rho(\boldsymbol\xi)
 =\sqrt{\left(\frac{s_1(1+t)+s_2(1-t)}{2}\right)^2
                   -s_1s_2(p_2^2+p_3^2)}.
\end{equation}
Thus $\widetilde\rho(\boldsymbol\xi)=\rho(\boldsymbol\xi)$
for $\boldsymbol\xi\in S^4$. To use the same formula on $D^5$,
we must first show that its square root remains real when
$\norm{\boldsymbol\xi}\leq1$.

\begin{lemma}\label{lem:second-radius}
Formula~\eqref{eq:rho-tilde} extends $\widetilde\rho$ from $S^4$
to a nonnegative real-valued function on $D^5$, and for every
$\boldsymbol\xi\in D^5$,
\begin{equation}\label{eq:radius-difference}
 \widetilde\rho(\boldsymbol\xi)^2-\rho(\boldsymbol\xi)^2
 =s_1s_2\bigl(1-\norm{\boldsymbol\xi}^2\bigr).
\end{equation}
\end{lemma}
\begin{proof}
We use $(a+b)^2-(a-b)^2=4ab$ to calculate
\[
 \begin{aligned}
 \widetilde\rho(\boldsymbol\xi)^2-\rho(\boldsymbol\xi)^2
 &=\left(\frac{s_1(1+t)+s_2(1-t)}{2}\right)^2
   -\left(\frac{s_1(1+t)-s_2(1-t)}{2}\right)^2\\
 &\quad-s_1s_2(p_0^2+p_1^2+p_2^2+p_3^2)\\
 &=s_1s_2(1-t^2-p_0^2-p_1^2-p_2^2-p_3^2)\\
 &=s_1s_2\bigl(1-\norm{\boldsymbol\xi}^2\bigr).
 \end{aligned}
\]
Since $s_1s_2\geq0$ and $\norm{\boldsymbol\xi}\leq1$ on $D^5$,
\[
 \widetilde\rho(\boldsymbol\xi)^2
 \geq\rho(\boldsymbol\xi)^2\geq0,
\]
so the square root defining $\widetilde\rho$ is real throughout
$D^5$.
\end{proof}

On $S^4$ we have $\widetilde\rho(\boldsymbol\xi)=\rho(\boldsymbol\xi)$,
so~\eqref{eq:Cplus} gives
\[
 \Cplus(A)=\{h(\boldsymbol\xi)+i\widetilde\rho(\boldsymbol\xi):
                         \boldsymbol\xi\in S^4\}.
\]
To pass from this sphere representation to the ball, we apply
the same composition argument as in
Proposition~\ref{thm:ball-representations}, now with the affine map
$\boldsymbol\xi\mapsto\bigl(h(\boldsymbol\xi),t,p_2,p_3\bigr)$.
In this case, $\widetilde\rho$ is independent of $p_0,p_1$;
with $t,p_2,p_3$ fixed, keeping $h$ unchanged again leaves a
nonzero direction in which these two remaining coordinates
can vary without changing the complex value.
It follows that $h+i\widetilde\rho$ has the same set of values
on $D^5$ as on $S^4$, yielding the following result.

\begin{corollary}\label{thm:ball-second-radius}
For the fixed matrix $A=H+Si$,
\begin{equation}\label{eq:ball-second-radius}
 \Cplus(A)=\{h(\boldsymbol\xi)+i\widetilde\rho(\boldsymbol\xi):
                                  \boldsymbol\xi\in D^5\}.
\end{equation}
\end{corollary}

We now examine how the two radius functions behave under convex
combinations in $D^5$. 

\begin{lemma}\label{lem:radius-convexity}
The function $\rho$ is convex and the function $\widetilde\rho$
is concave on $D^5$.
\end{lemma}
\begin{proof}
Let $\boldsymbol\xi_j=(t_j,\mathbf p_j)\in D^5$ for $j=1,2$,
where $\mathbf p_j=(p_{0,j},p_{1,j},p_{2,j},p_{3,j})$.
For $0\leq\theta\leq1$, their convex combination is
\[
 \boldsymbol\xi_\theta
 =\theta\boldsymbol\xi_1+(1-\theta)\boldsymbol\xi_2
 =(t_\theta,\mathbf p_\theta)\in D^5.
\]
Its coordinates satisfy $t_\theta=\theta t_1+(1-\theta)t_2$
and $p_{k,\theta}=\theta p_{k,1}+(1-\theta)p_{k,2}$ for
$k=0,1,2,3$.

The formula for $\rho$ can be written as
\[
 \rho(\boldsymbol\xi)=
 \left\lVert\left(
 \frac{s_1(1+t)-s_2(1-t)}{2},\,
 \sqrt{s_1s_2}\,p_0,\,\sqrt{s_1s_2}\,p_1
 \right)\right\rVert.
\]
Each of the three entries inside the norm is affine in
$\boldsymbol\xi$, so its value at $\boldsymbol\xi_\theta$ is
the corresponding convex combination of its values at
$\boldsymbol\xi_1$ and $\boldsymbol\xi_2$. Applying the triangle
inequality therefore gives
\[
 \rho(\boldsymbol\xi_\theta)
 \leq\theta\rho(\boldsymbol\xi_1)
       +(1-\theta)\rho(\boldsymbol\xi_2),
\]
which proves convexity.

To prove concavity of $\widetilde\rho$, put
\[
 y_j=\widetilde\rho(\boldsymbol\xi_j)\quad\text{for }j=1,2,
 \qquad y_\theta=\theta y_1+(1-\theta)y_2.
\]
We need to show that
$y_\theta\leq\widetilde\rho(\boldsymbol\xi_\theta)$.
By~\eqref{eq:rho-tilde},
\[
 \left\lVert\left(
 y_j,\,\sqrt{s_1s_2}\,p_{2,j},\,\sqrt{s_1s_2}\,p_{3,j}
 \right)\right\rVert
 =\frac{s_1(1+t_j)+s_2(1-t_j)}{2},
\]
where the right-hand side is nonnegative because
$-1\leq t_j\leq1$ and $s_1,s_2\geq0$.
Applying the triangle inequality to the convex combination of
these two vectors gives
\[
 \begin{aligned}
 &\left\lVert\left(
 y_\theta,\,\sqrt{s_1s_2}\,p_{2,\theta},\,
 \sqrt{s_1s_2}\,p_{3,\theta}\right)\right\rVert\\
 &\quad\leq
 \theta\left\lVert\left(
 y_1,\,\sqrt{s_1s_2}\,p_{2,1},\,\sqrt{s_1s_2}\,p_{3,1}
 \right)\right\rVert
 +(1-\theta)\left\lVert\left(
 y_2,\,\sqrt{s_1s_2}\,p_{2,2},\,\sqrt{s_1s_2}\,p_{3,2}
 \right)\right\rVert\\
 &\quad=
 \theta\frac{s_1(1+t_1)+s_2(1-t_1)}{2}
 +(1-\theta)\frac{s_1(1+t_2)+s_2(1-t_2)}{2}\\
 &\quad=\frac{s_1(1+t_\theta)+s_2(1-t_\theta)}{2}.
 \end{aligned}
\]
Squaring and rearranging yields
\[
 y_\theta^2
 \leq\left(\frac{s_1(1+t_\theta)+s_2(1-t_\theta)}{2}\right)^2
      -s_1s_2(p_{2,\theta}^2+p_{3,\theta}^2)
 =\widetilde\rho(\boldsymbol\xi_\theta)^2.
\]
Since both quantities are nonnegative, we obtain
\[
 \theta\widetilde\rho(\boldsymbol\xi_1)
 +(1-\theta)\widetilde\rho(\boldsymbol\xi_2)
 =y_\theta\leq\widetilde\rho(\boldsymbol\xi_\theta),
\]
which proves concavity.
\end{proof}

We now consider the parameters that give the same real part
in $\Cplus(A)$. For $a\in h(D^5)$, define the
\emph{fibre}
\[
 D_a=\{\boldsymbol\xi\in D^5:h(\boldsymbol\xi)=a\}.
\]
The values of $\rho$ on this set describe the imaginary parts
of the points of $\Cplus(A)$ with real part $a$. We show that
these values fill an interval.

\begin{lemma}\label{lem:fibres}
For every $a\in h(D^5)$, the fibre $D_a$ is nonempty, compact
and convex, and
\[
 \rho(D_a)=\widetilde\rho(D_a)
 =\{y\in\RR:a+iy\in\Cplus(A)\}
\]
is a nonempty compact interval.
\end{lemma}
\begin{proof}
The condition $a\in h(D^5)$ makes $D_a$ nonempty. Since $h$
is continuous and affine, its level set $h^{-1}(\{a\})$ is
closed and convex, so $D_a=D^5\cap h^{-1}(\{a\})$ is compact
and convex. In particular, $D_a$ is connected, and continuity
of $\rho$ implies that $\rho(D_a)$ is a nonempty compact
connected subset of $\RR$, hence a compact interval.
Finally, restricting~\eqref{eq:ball-representations} to points
with real part $a$ gives $\rho(D_a)$ as their set of imaginary
parts, while~\eqref{eq:ball-second-radius} gives
$\widetilde\rho(D_a)$ as that same set.
\end{proof}

To determine how the endpoints of these intervals vary with
$a$, we use the following general fact about extrema over
fibres of an affine map, a form of partial minimization;
see \cite[Section~3.2.5]{BV}.

\begin{lemma}\label{lem:fibre-extrema}
Let $K\subseteq\RR^d$ be nonempty, compact and convex, and let
$h:\RR^d\to\RR$ be affine. If $f,g:K\to\RR$ are continuous,
with $f$ convex and $g$ concave, then the functions
\[
 m(a)=\min_{\substack{\boldsymbol\xi\in K\\h(\boldsymbol\xi)=a}}
                  f(\boldsymbol\xi),\qquad
 M(a)=\max_{\substack{\boldsymbol\xi\in K\\h(\boldsymbol\xi)=a}}
                  g(\boldsymbol\xi)
\]
are convex and concave, respectively, on $h(K)$.
\end{lemma}
\begin{proof}
For $a\in h(K)$, the fibre $K\cap h^{-1}(\{a\})$ is nonempty
and compact, so both extrema are attained. Given $a_1,a_2\in h(K)$,
choose minimizers $\boldsymbol\xi_1,\boldsymbol\xi_2$ with
$h(\boldsymbol\xi_j)=a_j$ and $f(\boldsymbol\xi_j)=m(a_j)$.
For $0\leq\theta\leq1$, convexity of $K$ and affineness of $h$
place $\theta\boldsymbol\xi_1+(1-\theta)\boldsymbol\xi_2$
in the fibre over $a_\theta=\theta a_1+(1-\theta)a_2$. Hence
\[
 \begin{aligned}
 m(a_\theta)
 &\leq f\bigl(\theta\boldsymbol\xi_1+(1-\theta)\boldsymbol\xi_2\bigr)\\
 &\leq\theta f(\boldsymbol\xi_1)+(1-\theta)f(\boldsymbol\xi_2)
 =\theta m(a_1)+(1-\theta)m(a_2).
 \end{aligned}
\]
Thus $m$ is convex. Applying the same argument to the convex
function $-g$ shows that $-M$ is convex, so $M$ is concave.
\end{proof}

By Lemma~\ref{lem:fibres}, the two radius functions have the same
set of values on each fibre $D_a$, so either function gives the
same minimum and the same maximum there. We choose $\rho$ to
describe the minimum because it is convex, and $\widetilde\rho$
to describe the maximum because it is concave. These choices
allow us to apply Lemma~\ref{lem:fibre-extrema} to obtain a convex
lower boundary function and a concave upper boundary function,
as stated in the following proposition.

\begin{proposition}\label{prop:boundary-representation}
For $a\in h(D^5)$, define
\[
 m(a)=\min_{\boldsymbol\xi\in D_a}\rho(\boldsymbol\xi),
 \qquad
 M(a)=\max_{\boldsymbol\xi\in D_a}\widetilde\rho(\boldsymbol\xi).
\]
Then $m$ is convex and $M$ is concave on $h(D^5)$, and
\begin{equation}\label{eq:boundary-representation}
 \Cplus(A)=\{a+iy:a\in h(D^5),\ m(a)\leq y\leq M(a)\}.
\end{equation}
\end{proposition}
\begin{proof}
Both radius functions are continuous on $D^5$ by their formulas,
and Lemma~\ref{lem:radius-convexity} gives their convexity and
concavity. Applying Lemma~\ref{lem:fibre-extrema} with $K=D^5$,
$f=\rho$ and $g=\widetilde\rho$, we conclude that $m$ is convex
and $M$ is concave. Moreover, Lemma~\ref{lem:fibres} shows that
$\rho(D_a)=\widetilde\rho(D_a)$ is an interval, whose minimum
is $m(a)$ and whose maximum is $M(a)$. This gives the stated
representation of $\Cplus(A)$.
\end{proof}

These lower and upper boundaries give the final convexity argument.

\begin{proof}[Proof of Theorem~\ref{thm:upper-bild-convexity}]
We first consider $A=H+Si\in M_2(\HH)$, where $H=H^*$ and
$S=\diag(s_1,s_2)$ with $s_1,s_2\geq0$.
Take $z_j=a_j+ib_j\in\Cplus(A)$ for $j=1,2$.
The representation~\eqref{eq:boundary-representation} means that
$a_j\in h(D^5)$ and $m(a_j)\leq b_j\leq M(a_j)$ for each $j$.
For $0\leq\theta\leq1$, set
$a_\theta=\theta a_1+(1-\theta)a_2$ and
$b_\theta=\theta b_1+(1-\theta)b_2$.
The set $h(D^5)$ is an interval because $h$ is affine and $D^5$
is convex, so $a_\theta\in h(D^5)$.
The convexity of $m$, the concavity of $M$ and the bounds on
$b_1,b_2$ give
\[
 \begin{aligned}
 m(a_\theta)
 &\leq\theta m(a_1)+(1-\theta)m(a_2)\leq b_\theta\\
 &\leq\theta M(a_1)+(1-\theta)M(a_2)\leq M(a_\theta).
 \end{aligned}
\]
These inequalities show that $m(a_\theta)\leq b_\theta\leq M(a_\theta)$.
Together with $a_\theta\in h(D^5)$, this is precisely the condition
in~\eqref{eq:boundary-representation} for $a_\theta+ib_\theta$
to belong to $\Cplus(A)$. Hence
\[
 \theta z_1+(1-\theta)z_2
 =a_\theta+ib_\theta\in\Cplus(A).
\]
This proves that $\Cplus(A)$ is convex.

Proposition~\ref{thm:coordinate-model} identifies
$\Bplus(A)=\Cplus(A)$ for matrices in the chosen form.
Since every $2\times2$ quaternionic matrix can be brought to
this form by a unitary similarity that preserves the upper bild,
convexity holds for every $2\times2$ matrix.
Lemma~\ref{lem:compression} then gives the result in every
dimension.
\end{proof}


\begin{thebibliography}{GR}
\raggedright

\bibitem[BV]{BV}
S. Boyd and L. Vandenberghe,
\emph{Convex Optimization}, Cambridge University Press, 2004.
\href{https://web.stanford.edu/~boyd/cvxbook/bv_cvxbook.pdf}{Online edition}.

\bibitem[GR]{GR}
K. E. Gustafson and D. K. M. Rao,
\emph{Numerical Range: The Field of Values of Linear Operators and Matrices},
Universitext, Springer-Verlag, New York, 1997.
\href{https://doi.org/10.1007/978-1-4613-8498-4}{doi:10.1007/978-1-4613-8498-4}.

\bibitem[ST]{ST}
W. So and R. C. Thompson,
\emph{Convexity of the upper complex plane part of the numerical range
of a quaternionic matrix}, Linear Multilinear Algebra \textbf{41}(4)
(1996), 303--365.
\href{https://doi.org/10.1080/03081089608818484}{doi:10.1080/03081089608818484}.

\bibitem[Zh]{Zh}
F. Zhang,
\emph{Quaternions and matrices of quaternions},
Linear Algebra Appl. \textbf{251} (1997), 21--57.
\href{https://doi.org/10.1016/0024-3795(95)00543-9}{doi:10.1016/0024-3795(95)00543-9}.

\end{thebibliography}
\end{document}